\documentclass[a4paper,12pt]{article}
\usepackage{amssymb}
\usepackage{amsmath}
\usepackage{graphicx}
\usepackage{epstopdf}

\newcommand{\eneqa}{\end{eqnarray}}
\newcommand{\begeqaet}{\begin{eqnarray*}}
\newcommand{\eneqaet}{\end{eqnarray*}}
\newcommand{\be}{\begin{equation}}
\newcommand{\ee}{\end{equation}}

\newcommand{\rn}{\rbig^N}

\newcommand{\rbig}{{\mathbb{R}}}

\newcommand{\ed}{\end{document}}

\newcommand{\befr}{\begin{frame} }
\newcommand{\efr}{\end{frame}}

\def\beq{\begin{equation}}
\def\eeq{\end{equation}}

\newtheorem{theo}{Theorem}

\newtheorem{prop}{Proposition}[section]

\newcommand{\begeqa}{\begin{eqnarray}}

\begin{document}

\begin{center}{\bf\large Uniform bounds via regularity estimates for elliptic PDE with critical growth in the gradient}
\end{center}

 \begin{center}
Boyan SIRAKOV\footnote{e-mail : bsirakov@mat.puc-rio.fr}\\
PUC-Rio, Departamento de Matematica,\\ Rio de Janeiro - CEP 22451-900, BRASIL \\
\end{center}

{\small \noindent{\bf Abstract}. We prove  non-uniqueness and study the behaviour of viscosity solutions of a class of uniformly
elliptic fully nonlinear equations of Hamilton-Jacobi-Bellman-Isaacs type, with  quadratic growth in the gradient. The crucial a priori bound for the solutions is proved through an argument which uses a boundary growth lemma, and consequences such as boundary "half"-Harnack inequalities, which are of independent interest. Our results are new even for linear equations.}
%

\section{Introduction and Main Results}

This note is an account of some recent results on solvability and multiplicity of solutions of the Dirichlet problem for a uniformly elliptic operator in which the first order term has the same scaling with respect to dilations as the second order term. This work can be considered as a continuation of the paper \cite{S1}, and is also strongly motivated by a number of very recent papers on a particular equation of our type, \cite{JS}, \cite{ACJT}, \cite{CJ}. The purpose of this note is to state our main theorems in a particular but typical case, and give a sketch of their proofs; the full presentation will appear in the forthcoming work \cite{S2}.

Apart from their appearance in a number of applications, this type of PDE is of theoretical importance, since it represents a class of equations which is {\it invariant} with respect to diffeomorphic changes of dependent and independent variable. An additional difficulty in the study of these equations is the "critical" behaviour of the gradient, in the sense that the first-order term has the same scaling with respect to dilations as the second-order term, and thus does not scale out after a "zoom" in a point.

Even though the results we obtain are valid for general Isaacs operators as in \cite{S1}, here for simplicity we will state them only for the equation

\be\label{main}
\left\{
\begin{array}{rclcc}
-L_0 u  &=& c(x)u\: + <\!M(x)\nabla u, \nabla u\!> +\:h(x)&\mbox{ in }& \Omega\\
u&=&g(x)&\mbox{ on }&\partial \Omega
\end{array}\right.
\ee
where $\Omega$ is a $C^2$-smooth bounded domain in $\rn$, and
$$
L_0u = a_{ij}(x)\partial _{ij}u + b_i(x)\partial_i u
$$
is a linear uniformly elliptic  operator whose coefficients have the necessary regularity to ensure satisfactory theory for the linear problem ($M=0$). For simplicity here we will assume all coefficients in \eqref{main} are bounded functions,  $A$ is continuous and $\lambda I\le A(x)\le \Lambda I$, $|M(x)|\le M_0 $, $0<\lambda\le \Lambda$, $M_0\ge0$.

The equation \eqref{main} can be thought of as the "closure" of the class of linear elliptic equations with respect to diffeomorphic changes of the dependent variable ($u\to v=\psi(u)$) and of the independent variable ($x\to y=\Psi(x)$).

  We have separated the zero order term in the elliptic operator, since the sign and size of $c(x)$ matter, both for the solvability and for the uniqueness of solutions.

When reduced to \eqref{main}, the main theorem in \cite{S1} states the following.
\begin{enumerate} \item The problem \eqref{main} has a unique solution if $c(x)\le -c_0<0$.
\item There exists $\delta_0>0$ depending on $\lambda$, $\Lambda$, $\|b\|_n$, diam$(\Omega)$, such that if
$$
\max\:\{\, M_0\|h\|_{L^n},\: \|c^+\|_{L^n},\: M_0\left(\max_{\partial\Omega} g\right)\|c^+ \|_{L^n}\:\}\le \delta_0
$$
then \eqref{main} has a solution. This hypothesis cannot be improved in general. The solution is unique if $c^+=0$.
\item If $c\equiv c_0\in (0,\delta_0)$, $L_0=\Delta$, $M=\mu_0 I$, $\mu_0>0$, $g=h=0$, then \eqref{main} has at least two solutions.
\end{enumerate}

The possibility of extending the result in 3. above, that is, to show non-uniqueness in general when $c$ is not nonpositive was left as an open question in \cite{S1}. This question was taken up in a number of recent works, which we now describe in a little more detail. In the following we also assume for simplicity that the boundary data $g=0$ on $\partial \Omega$, and that $c(x)\gneqq0$ in $\Omega$. Of course we always need to assume that the problem with $c\equiv0$ has a solution (as we know, such a solution is unique, and its existence is implied by an upper bound on $M_0\|h\|_{L^n}$).

In \cite{CJ} it was shown that  if the second order operator $L_0$ is the Laplacian, and $M(x) = \mu(x)I$, where the function $\mu$ is such that $|\mu(x)|\ge \mu_1>0$, then \eqref{main} has at least two solutions if $c(x)\le C_0$, for some explicit constant $C_0$. That paper extends the earlier work \cite{JS} and further develops the topological degree method used in \cite{S2} and \cite{ACJT}. The behaviour of the solutions with respect to a parameter which measures the size of $c(x)$ is also studied in \cite{CJ}, through a Rabinowitz-type bifurcation method.

The following diagrams from \cite{CJ} are worth reproducing. Consider the model problem for a given parameter $\lambda\in \mathbb{R}$
\be\label{model}\left\{
\begin{array}{rclcc}
-\Delta u  &=&\lambda c(x)u\: + \mu(x)|\nabla u|^2 +\:h(x)&\mbox{ in }& \Omega\\
u&=&0&\mbox{ on }&\partial \Omega
\end{array}\right.
\ee

\begin{figure}[h]
\includegraphics[width=4.5cm, height=4cm]{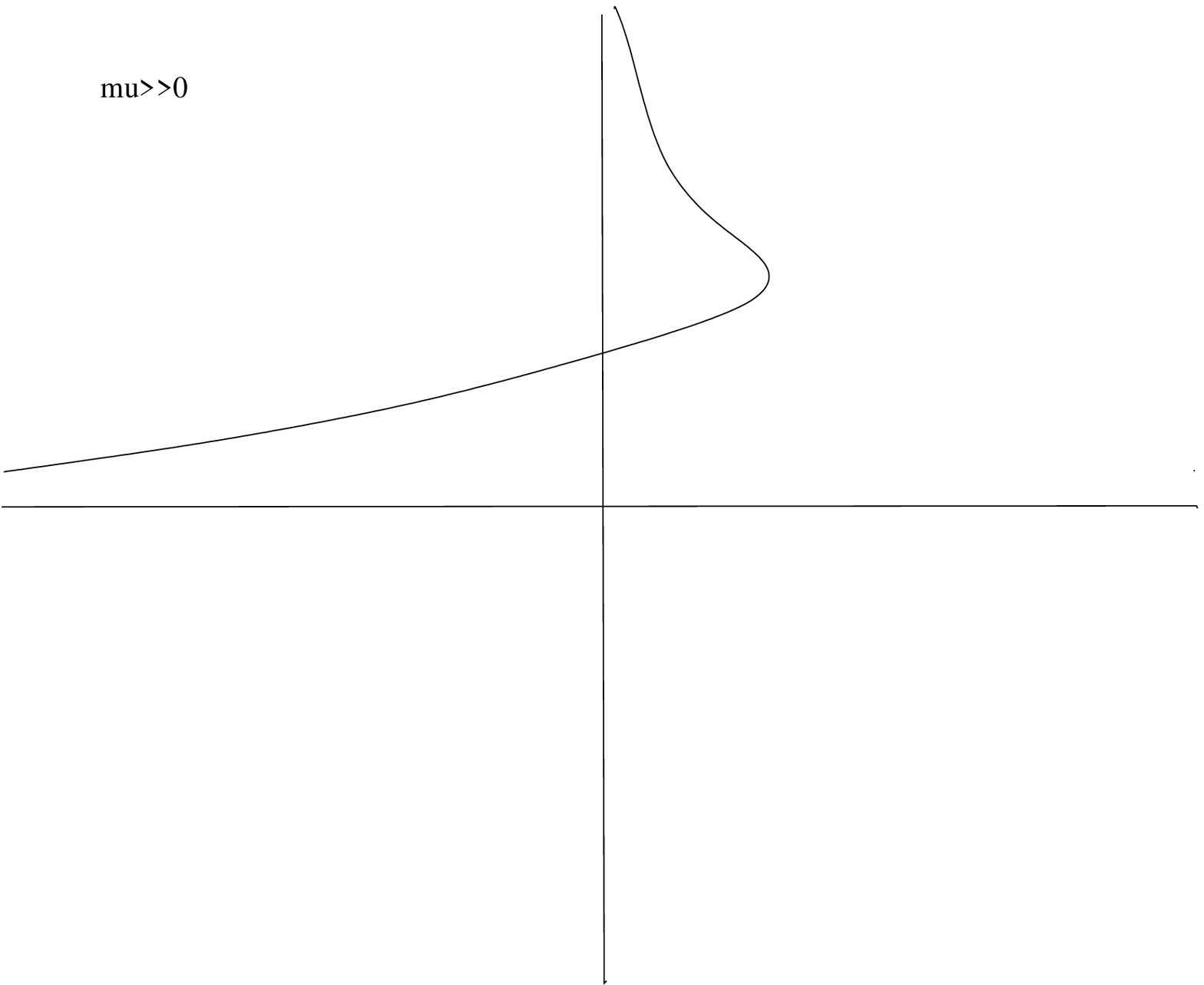}
\includegraphics[width=4.5cm, height=4cm]{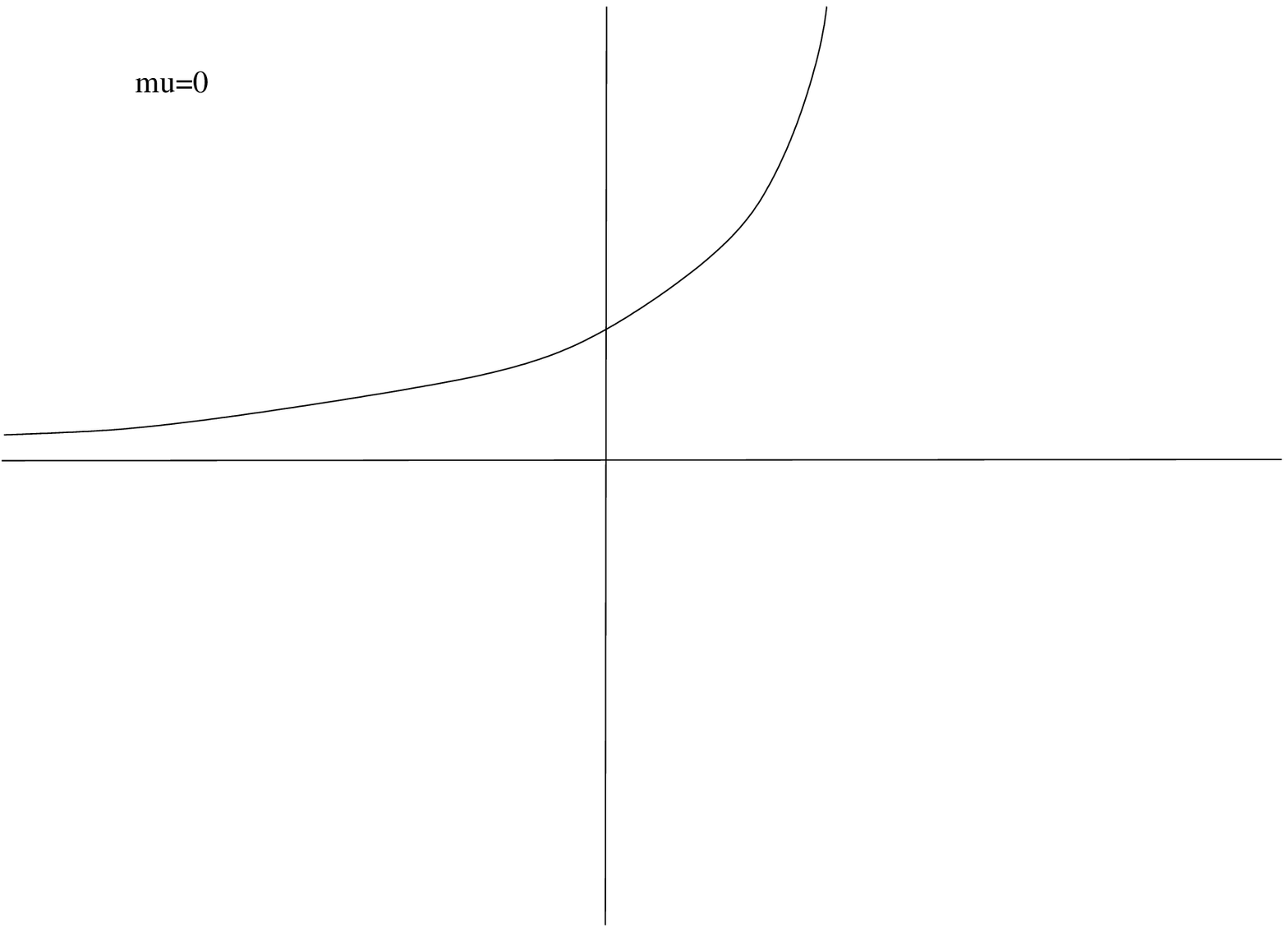}
\includegraphics[width=4.5cm, height=4cm]{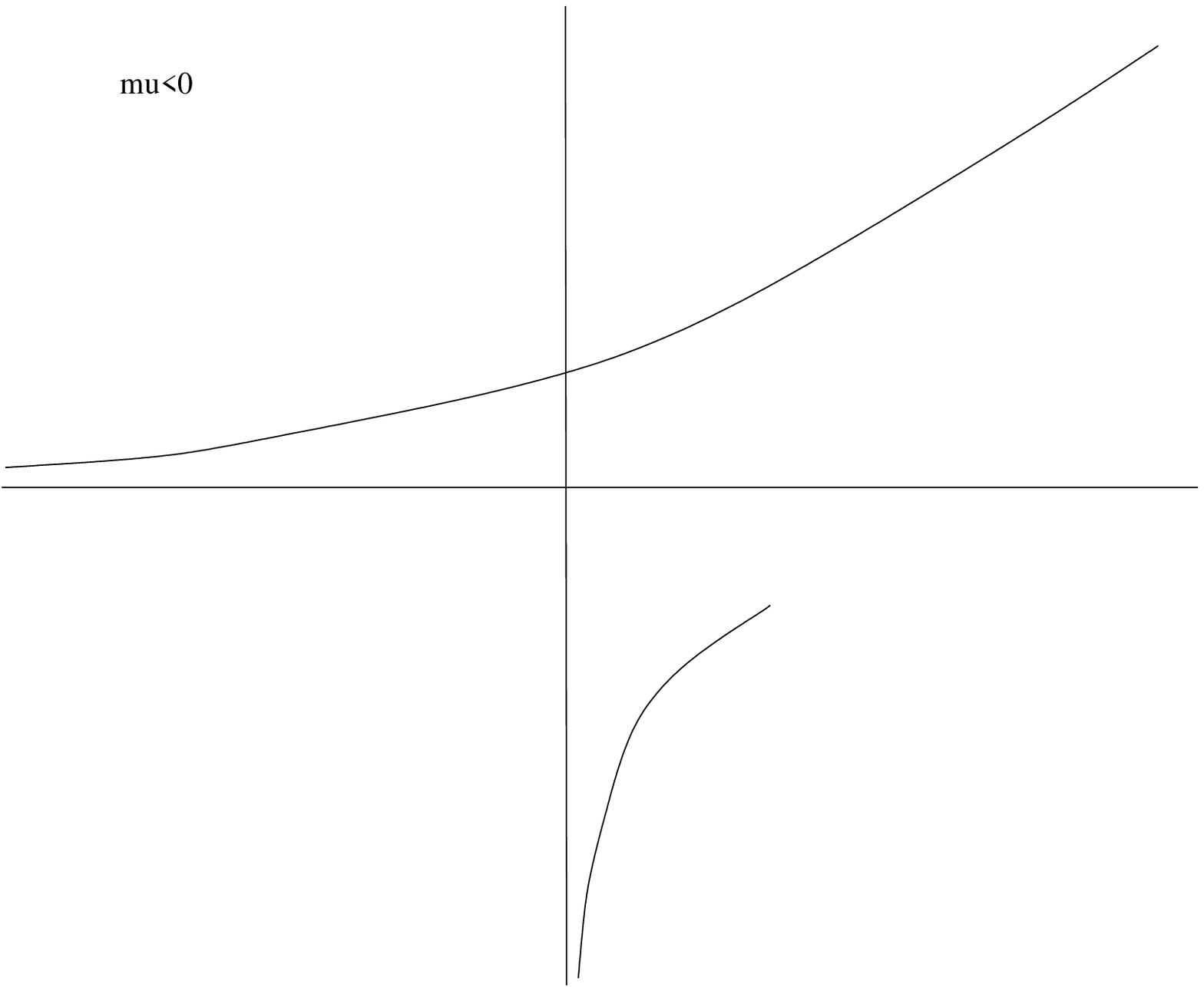}
\label{fig2}
\end{figure}

Here we have mapped the parameter $\lambda$ on the horizontal axis, and the value of the solution of \eqref{model} for some point in the domain on the vertical axis.
These diagrams correspond to the cases $\mu>>0$, $\mu=0$ and $\mu<<0$ respectively, and are valid for instance if $h\gneqq0$. The main difficulty in establishing their validity is to obtain uniform a priori bounds for the solutions, for all $\lambda$ in any interval of the form $[0,A]$ if $\mu<0$, and  in any interval of the form $[a,A]$, $0<a\le A<\infty$, if $\mu>0$.

As can be seen from the above pictures, it might be complicated to understand the behaviour of solutions if we assume only $\mu(x)\gneqq0$ (instead of $\mu(x)\ge\mu_1>0$). This question for the model equation \eqref{model} was studied in the very recent work \cite{So}, where it was shown that a necessary condition for an a priori bound when $\mu(x)\ge0$ is that the intersection of the supports of $c(x)$ and $\mu(x)$ contains a ball. It was also shown that this condition is sufficient if (i) the spatial dimension $n=2$; (ii) or if we assume $n=3$ and a growth assumption on $c$ and $\mu$ close to $\partial \Omega$; (iii) or if $n\le 4$ and $\mu(x)\ge \mu_1>0$ on supp$(c)$.

In all previous works the results were restricted to the model equation \eqref{model}, since the obtention of the a priori bounds crucially depended on the fact that the second order operator in the equation is the Laplacian. In \cite{JS} a variational argument of "mountain-pass" type was used, in \cite{ACJT} and \cite{CJ} the weak-Sobolev formulation of the system \eqref{feq}-\eqref{seq} below is tested by the first eigenfunction and the solution itself, and in \cite{So} interpolation and elliptic estimates for the Laplacian in weighted Lebesgue spaces were used.

It is our goal here to present another method for obtaining a priori bounds, which depends only on the uniform ellipticity of the operator, and thus gives the above described results for any operator $L_0$. We prove the upper bound by using a growth lemma (quantitative strong maximum principle) combined with a weak Harnack inequality and a local maximum principle. A supplementary difficulty is that boundary versions of these results do not appear to be available, and need to be proved. This can be done by a modification of the original Krylov's argument for proving boundary regularity for elliptic equations. Furthermore, when $\mu(x)\ge 0$ we prove the uniform bound under the hypothesis $\mu(x)\ge \mu_1>0$ on supp$(c)$ in any dimension, thus extending and joining the results for the Laplacian from \cite{CJ} and \cite{So} (see (iii) above).

It is also interesting that the  bound in the case $\mu\le 0$ turns out to be equivalent to the validity of the strong maximum principle for a coercive elliptic inequality with a superlinear nonlinearity which satisfies the Vazquez condition.

\begin{theo}
Assume $M(x)\ge \mu_1I$ or $M(x)\le -\mu_1I$ on $\mathrm{supp}(c)$, $\mu_1>0$.

 If the problem \eqref{main} with $c=0$ has a solution $u_0$, then \eqref{main} has at least two solutions for $0\lneqq c(x)\le \delta_0$, where $\delta_0$ depends on the ellipticity, the upper bounds for the coefficients of the equation, and on the domain.

 If $M(x)\ge \mu_1I$  on $\mathrm{supp}(c)$ and $u_0\gneqq0$ then the first diagram in Fig.1 is valid.
 If $M(x)\le -\mu_1I$  on $\Omega$ and $u_0\gneqq0$ then the third diagram is valid.
\end{theo}

In the proof of this theorem we use the following boundary estimates which are clearly of independent interest. We assume we have a domain with a flat portion of the boundary $\{x_n=0\}$, and denote with $B_R^+$ a half-ball in the domain centered at a point on this portion. We set $B_R^0=\partial B_R^+\cap\{x_n=0\}$.

\begin{prop} (boundary quantitative strong maximum principle, BQSMP) Assume $u$ is a solution of $-L_0u\ge 0$, $u\ge0$ in $B_2^+$, and $u=0$ on $B_2^0$. Then there exists $\varepsilon>0$ depending on $\lambda$, $\Lambda$ and $\|b\|_{L^p}$ such that

$$
\inf_{B_1^+} \frac{u}{x_n} \ge c\left( \int_{B_{1}^+} (-L_0u)^{\varepsilon}\right)^{1/{\varepsilon}}.
$$
\end{prop}

\begin{prop} (boundary weak Harnack inequality, BWHI) Assume that $-L_0u\ge f$, $u\ge 0$ in $B_2^+$, $u=0$ on $B_2^0$. Then there exists $\varepsilon>0$ such that
$$
 \left( \int_{B_{3/2}^+} \left(\frac{u}{x_n}\right)^\varepsilon\right)^{1/\varepsilon}\le C \left(\inf_{B_{3/2}^+} \frac{u}{x_n} + \|f^-\|_{L^n(B_2^+)}\right).
$$
\end{prop}

\begin{prop} (boundary local maximum principle, BLMP) Assume that $-L_0u\le d(x)u + f$ in $B_2^+$, $u=0$ on $B_2^0$, $d\in L^{q}$, for some $q>n$. Then for each $p>0$
$$
\sup_{B_{1}^+} \frac{u^+}{x_n}\le C\left(\left( \int_{B_{3/2}^+} (u^+)^p\right)^{1/p} + \|f^+\|_{L^n(B_2^+)}\right)
$$
\end{prop}

 The first two of these propositions can be proved by using Krylov's original idea for proving boundary regularity, which consists in writing a equation for $u/x_n$ that degenerates in a special way, and for which the fundamental "growth lemma" can still be proved. See \cite{S2} for details. An interior QSMP in this convenient form can be found in \cite{K}, the interior versions of WHI and LMP for instance in \cite{GT}. Of course the BLMP above is a  consequence of the boundary Lipschitz estimate and the usual LMP ($u^+$ is a subsolution in a full ball if we extend $u=0$ outside $\Omega$).

It is interesting to notice that the combination of BWHI and BLMP yields the full boundary Harnack inequality, a result which is very well known, and goes back to \cite{CCCC}, \cite{B}. However, in all texts where this inequality has appeared it has been proved by a method different from the above splitting into separate results for supersolutions and subsolutions.

It is also worth noting that inequalities  called "boundary weak Harnack inequalities" have appeared in a different form, which is sufficient for a proof of H\"older regularity up to the boundary (see for instance \cite{GT}, Theorems 8.26 and 9.27). However, that inequality is void for a function which vanishes on the boundary and does not imply the full boundary Harnack inequality.
\bigskip

\noindent{\it Sketch of the proof of the theorem}. First, we observe that if $c<\lambda_1(L_0,\Omega)$ then we can assume without loss of generality that $M(x)\ge0$, $h\ge0$, $u_0\ge0$ and we can search for positive solutions of our equation. This is done by observing that the difference of $u$ and the solution of $L_0 v = c(x) v + h$ satisfies an equation (with a modified operator $\tilde{L}_0$) in which these hypotheses are verified.

Next, we observe that if there exists an unbounded sequence of solutions, then there also exists a sequence of solutions $u_k$ and points $x_k$ such that $u_k(x_k)\to \infty$, and $x_k\to x_0\in\overline{\Omega}$ is such that for any ball $B$ centered at $x_0$ we have $\int_{B\cap\Omega} c>0$. This is because for any domain $G\subset\{c\equiv0\}$ we can apply the comparison principle in $G$ to the functions $u-\sup_{\partial G} u$ and $u_0-\inf_\Omega u_0$.

So we see we need to prove an a priori bound in a ball $B$ (resp. a half-ball if $x_0\in \partial \Omega$) such that $c\gneqq0$ in $B$ and $0<\mu_1\le \mu(x)\le \mu_2$ in $\Omega$, by the hypothesis of the theorem.

We make a classical exponential change of the dependent variable $u$, setting $v_i= \mu_i^{-1}(e^{\mu_i u}-1)$, $i=1,2$.

In $B$, for large values of $u$
\begeqa
-L_1 v_1&\ge &f_1(x,v_1)\sim c_0\,c(x)\,v_1\,\log(v_1)\label{feq}\\
-L_2 v_2&\le &f_2(x,v_2)\sim C_0\,c(x)\,v_2\,\log(v_2)\label{seq}\\
v_2&\sim&v_1^A,\qquad A=\mu_2/\mu_1\nonumber\\
v_i&=&0 \mbox{ on }\partial \Omega\cap B.\nonumber
\eneqa
To simplify, and since we only want to give the main idea of the proof here, we do not write the full expression of $f_i$; what actually matters is the above asymptotic behaviour for large values of $u$ (resp. $v_i$). We observe at this point that a priori bounds for this system of inequalities do not seem to be provable by classical methods such as the Gidas-Spruck rescaling method. This is because we may have $c(x_0)=0$, and, on the other hand, if we write this system as a system in only one of the functions $v_i$, then one of the operators will become degenerate.

Let us assume we are in the more difficult case, when $B$ is a half-ball, with $\partial \Omega\cap B\subset\{x_n=0\}$ (after a "flattening" diffeomorphic change of the independent variable, which only changes the coefficients of the equation).

Then the first inequality \eqref{feq} and the BQSMP imply
$$
\inf_{B_1^+} \frac{v_1}{x_n} \ge c_0\; \inf_{B_1^+} \frac{v_1}{x_n}\; \log\!\left( \inf_{B_1^+} \frac{v_1}{x_n}\right)\;
\left(\int_{B_1^+} c^\varepsilon x_n^{1+\varepsilon}\right)
$$

This implies that $\displaystyle \inf_{B_1^+} \frac{v_1}{x_n}\le C$. Then $-L_1 v_1\ge 0$ (which follows from \eqref{feq}) and the BWHI imply
$$
 \left( \int_{B_1^+} \left(v_1\right)^\varepsilon\right)^{1/\varepsilon}\le C \left( \int_{B_1^+} \left(\frac{v_1}{x_n}\right)^\varepsilon\right)^{1/\varepsilon}\le C.
$$

Now we apply the BLMP to the second inequality \eqref{seq}, in which $d(x)=C_0 c(x)\log(v_2)$ (note $\log(v_2)\sim A\log(v_1)$). Since for large values of $v_1$ we have $(\log(v_1))^{n+1}\le C v_1^\varepsilon$, we obtain from LMP and what we already proved
$$
\sup_{B_{1}^+} v_2\le \left( \int_{B_1^+} \left(v_2\right)^{\varepsilon/A}\right)^{A/\varepsilon} \le C \left( \int_{B_1^+} \left(v_1\right)^\varepsilon\right)^{A/\varepsilon} \le C.
$$
which establishes the required a priori bound, in terms of the ellipticity, bounds on the coefficients, and a lower bound on $\int_B c^\varepsilon$. We remark it would have been enough to use instead of the BQSMP a "growth lemma" which says that a positive supersolution is uniformly bounded below by a constant which depends on the measure of the set where the right-hand side of the inequality is strictly positive. 

Observe we did not need to reason by contradiction in order to obtain the a priori bound, nor had we to rescale the solutions or use any passage to the limit. We also do not need to assume any continuity of the coefficients of the equation, including the leading terms $a_{ij}$.

 It should be noted that {\it interior} Harnack inequalities (but not QSMP) have been used for proving a priori bounds on compact subsets of a domain, for instance in \cite{AL}, \cite{R}, in conjunction with  variational structure or rescaling arguments, which require more regularity and that an equation (as opposed to two one-sided inequalities) be satisfied by the solution.
\bigskip

It is worth observing that the results (i) and (ii) from the paper \cite{So} which we quoted above can be understood differently, by applying the above argument in which one uses the weak Harnack inequality applied to the function $u$ (which is itself superharmonic) instead of $v_1$, and the exact best exponent $\varepsilon$ in this inequality (which is known for the Laplacian), in order to see what is the best regularity one can get for $d(x) = c(x)u(x)$ in the inequality $-\Delta v_2\le f_2(x,v_2)\sim C c(x) u(x) v_2$, so that $d\in L^q$ with $q>n/2$. This leads us to the conjecture that (i) and (ii) from \cite{So}, as quoted above, cannot be improved.
\bigskip

Finally, let us sketch the proof of the a priori bound  in the "negative" case $M(x)\le -\mu_1 I $ in $\Omega$. In other words, we want to prove an upper bound for subsolutions of
\be\label{main2}
\left\{
\begin{array}{rclcc}
-L_0 u  &\le& c(x)u\: -\mu_1|\nabla u|^2 +\:h(x)&\mbox{ in }& \Omega\\
u&=&0&\mbox{ on }&\partial \Omega
\end{array}\right.
\ee

By setting $v= \mu_1^{-1}(1-e^{\mu_i u})$ we see that $v< 1$ is a subsolution of
\begeqaet
-L_1 v &\le& c (1-v) |\log(1-v)| + h(1-v) =:f(x,v)\mbox{ in } \Omega\\
v&=&0\mbox{ on }\partial \Omega
\eneqaet
Here the function $f(x,v(x))$ is {\it bounded} in $\Omega$ for $v_0\le v\le1$, where $v_0$ is any fixed subsolution of the above inequality. So by the standard Lipschitz bound at $\partial \Omega$ we have $v\le C\mathrm{dist}(x,\partial\Omega)$. Denote with $\bar{v}$ the supremum of all subsolutions $v\in (v_0,1)$ of the last inequality. Then by Perron's method $\bar{v}$ is also a subsolution. By the Lipschitz bound $v\not \equiv1$.

Now, the existence of an unbounded sequence of subsolutions of \eqref{main2} means precisely that $z=1-v\not \equiv0$ is a nonnegative supersolution of
$$
-L_2 z\ge c_0|\log(z)| z
$$
and $z$ vanishes somewhere in $\Omega$. This is impossible by Vazquez's strong maximum principle \cite{V} (and its extension to nonlinear inequalities, see for instance \cite{PS}, \cite{Q}), since one over the square root of the primitive of $t\log(t)$ is not integrable at zero.

\ed

{Some history}

\textsc{Kazdan---Kramer} 1975,
various partial results\bigskip

\textsc{Boccardo---Murat---Puel} 1980-1985\\
succeeded to obtain a full solvability result for the general class of divergence form operators, in the case when
$$
c(x)\le -c_0<0
$$\bigskip

\textsc{Ferone---Murat} 2000\\
solvability for the case
$$c\equiv0$$
and observation that only possible if $|Mh|$ small. \bigskip

\textsc{Barles-Murat, B.-Porretta}... 1995-2000 --- uniqueness. \bigskip

{\small
Related results by dall'Aglio, Giachetti and Puel; Maderna, Pagani and Salsa; Grenon, Murat and Porretta; Abdellaoui, dall'Aglio and Peral;  Abdel Hamid and Bidaut-V\'eron; Boccardo, Gallouet, Murat...}
\end{frame}

Results for the model equation with $L_0=\Delta$,\ \ \ $M=\mu(x) I$,

$$
-\Delta u  = c(x)u\: + \mu(x)|\nabla u|^2 +\:h(x)
$$\pause

\begin{block}

If $c<\lambda_1$ then $M\ge0$, $h\ge0$ and $u\ge 0$ is generic.\bigskip

$$
-\Delta u  = c(x)u\: + \mu(x)|\nabla u|^2 +\:h(x)
$$
$$
-\Delta \psi  = c(x)\psi +h(x)
$$\bigskip

Then $v=u-\psi$ is positive (MP) and satisfies
$$
-\Delta v +<\!\nabla\psi,\nabla v\!> = c(x)v\: + \mu(x)|\nabla v|^2 + \mu(x)|\nabla \psi|^2
$$\bigskip

modified elliptic operator.
\end{frame}

\begin{frame}{Some observations}

If lack of a priori bound, then lack of a priori bound on supp$(c)$.
\bigskip

In $\mathcal{O}$ open in $\{c\equiv0\}$, then in $\mathcal{O}$\\ both $u-\sup_{\partial\mathcal{O}} u$ and $u_0-\inf_{\Omega} u_0$ are solutions of

$$
-\Delta u  =  \mu(x)|\nabla u|^2 +\:h(x)
$$

and they compare on the boundary, so by MP

$$
\sup_{\mathcal{O}} u \le \sup_{\partial\mathcal{O}} u +2\|u_0\|_\infty
$$
\end{frame}

\begin{frame}{Some observations}

So if lack of a priori bound, then for some $x_0\in\overline{\Omega}$ and some (half-)ball $B$ around $x_0$
we have an unbounded in $B$ sequence of solutions of

$$
-L_0 u  = c(x)u\: + \mu(x)|\nabla u|^2 +\:h(x)
$$

and in $B$

$$
0<\mu_1\le \mu(x)\le \mu_2, \qquad c\gneqq0.
$$\pause

standard exponential change $v_i= \mu_i^{-1}(e^{\mu_i u}-1)$

In $B$, for large $u$
\begeqaet
-L_1 v_1&\ge &f_1(x,v_1)\sim c_0\,c(x)\,v_1\,\log(v_1)\\
-L_2 v_2&\le &f_2(x,v_2)\sim C_0\,c(x)\,v_2\,\log(v_2)\\
v_2&\sim&v_1^A
\eneqaet
\end{frame}

\begin{frame}

In $B$, for large $u$
\begeqaet
-L_1 v_1&\ge & c_0\,c(x)\,v_1\,\log(v_1)+h_1\\
-L_2 v_2&\le & C_0\,c(x)\,v_2\,\log(v_1)+h_2\\
v_2&\sim&v_1^A
\eneqaet

\begin{theorem} There exists a constant $C$ depending on $\lambda$, $\Lambda$, $\|b\|_p$, $\|c\|_\infty$, and a lower bound on $\int_B c$, such that $\|u\|_\infty\le C$.
\end{theorem}
\end{frame}

\begin{frame}

1. (QSMP) $-Lu\ge 0$, $u>0$ in $\Omega$ $\Longrightarrow$ for each compact $K\subset\Omega$
$$
\inf_K u \ge c\left( \int_K (-Lu)^{\varepsilon_0}\right)^{1/{\varepsilon_0}}.
$$\bigskip

2. (WHI) $-Lu\ge f$, $u>0$ in $\Omega$ $\Longrightarrow$ for each  $K\subset K^\prime\subset\Omega$
$$
 \left( \int_{K^\prime} u^\varepsilon\right)^{1/\varepsilon}\le C (\inf_K u + \|f\|_n).
$$\bigskip

3. (LMP) $-Lu\le f$ in $\Omega$ $\Longrightarrow$ for each  $K\subset K^\prime\subset\Omega$ and  $p>0$
$$
\sup_K u\le C\left(\left( \int_{K^\prime} u^p\right)^{1/p} + \|f\|_n\right)
$$
\pause\bigskip

Need boundary versions of these !!
\end{frame}

\begin{frame}

\end{frame}

\begin{frame}

\begeqaet
-L_1 v_1&\ge & c_0\,c(x)\,v_1\,\log(v_1)+h_1\\
-L_2 v_2&\le & C_0\,c(x)\,v_2\,\log(v_1)+h_2\\
v_2&\sim&v_1^A
\eneqaet

\end{frame}
\ed
\begin{thebibliography}{99}\small

\bibitem{AL} H. Amann, J. L\'opez-G\'omez, A priori bounds and multiple solutions for superlinear indefinite elliptic problems. J. Diff. Eq. 146 (2) (1998), 336-374.

\bibitem{ACJT} D. Arcoya, C. De Coster, L. Jeanjean, K. Tanaka, Continuum of solutions for an elliptic problem with critical growth in the gradient. J. Funct. Anal. 268 (8) (2015), 2298-2335.
    
\bibitem{B}   P.  Bauman,  Positive solutions of elliptic equations in nondivergence form and their adjoints. Ark. Mat. 22 (2) (1984),  153-173.
    

 \bibitem{BMP} L. Boccardo, F. Murat, J.-P. Puel,   Existence de solutions
faibles pour des \'equations elliptiques quasi-lin\'eaires \`a
croissance quadratique.  Res. Notes in Math. 84, Pitman, 1983,
19-73.
    
\bibitem{CCCC}   L. Caffarelli, E. Fabes, S. Mortola, S. Salsa,  Boundary behavior of nonnegative solutions of elliptic operators in divergence form. Indiana Univ. Math. J. 30 (4) (1981),  621-640.
    

 \bibitem{CJ} C. De Coster, L. Jeanjean,  Multiplicity results in the non-coercive case for an elliptic problem with critical growth in the gradient, preprint, arXiv:1507.04880.

\bibitem{Q} P. Felmer, A. Quaas, B. Sirakov,  Solvability of nonlinear elliptic equations with gradient terms. J. Diff. Eq. 254 (11) (2013), 4327-4346.

 
 \bibitem{GT} D. Gilbarg, N.S. Trudinger, { Elliptic Partial
Differential Equations of Second Order}, 2nd edition, Springer
Verlag.
 
\bibitem{JS} L. Jeanjean, B. Sirakov,  Existence and multiplicity for elliptic problems with quadratic growth in the gradient. Comm. Part. Diff. Eq. 38 (2013),  244-264.
    
\bibitem{K} N.V. Krylov,  Some Lp-estimates for elliptic and parabolic operators with measurable coefficients. Discr. Cont. Dyn. Syst. B 17 (6) (2012), 2073-2090.
   
   \bibitem{PS} { P. Pucci, J. Serrin}, {\ The maximum principle}. Progress in Nonlinear Differential Equations and their Applications, 73. Birkhauser Verlag, Basel, 2007.
    
 \bibitem{R}   D. Ruiz,  A priori estimates and existence of positive solutions for strongly nonlinear problems. J. Diff. Eq. 199 (1) (2004),  96-114.
    
\bibitem{S1} B.  Sirakov, Solvability of uniformly elliptic fully nonlinear PDE. Arch. Ration. Mech. Anal. 195 (2) (2010),  579-607.
  
 \bibitem{S2} B. Sirakov, Solvability of uniformly elliptic fully nonlinear PDE with critical growth in the gradient. In preparation.
  
\bibitem{So}  P. Souplet  A priori estimates and bifurcation of solutions for a noncoercive elliptic equation with critical growth in the gradient, preprint,
     arXiv:1411.0884.
     
\bibitem{V}     V\'azquez, J. L. A strong maximum principle for some quasilinear elliptic equations. Appl. Math. Optim. 12 (3) (1984),  191-202.
\end{thebibliography}
